\documentstyle[11pt,epsfig]{article}
\begin{document}
\bigskip

\begin{center}
{\Large \bf `Fair' Partitions of Polygons : An Introduction}\\
{\bf R Nandakumar$^\star$ and N Ramana Rao$^\dagger$}\\
\vskip .5cm
{\normalsize \it $\star$ R.K.M.Vivekananda University, Belur Math, Howrah, India - nandacumar@gmail.com}\\
{\normalsize \it $\dagger$ nallurirr@gmail.com}\\
\vskip 0.2cm
\vskip 1cm

{\bf ABSTRACT}
\end{center}
We address the question: Given a positive integer $N$, can any 2D convex polygonal region be partitioned into $N$ convex pieces such that all pieces have the same area and same perimeter?\\
The answer to this question is easily `yes' for $N$=2. We give a proof that the answer is `yes' for $N$=4 and also discuss higher powers of 2. \\

\section{Introduction - the Problem and a Conjecture}

A recently studied problem in Combinatorial Geometry is to partition a given convex polygon into convex pieces of equal area with the pieces also sharing the boundary of the input polygon equally (\cite{ref1} and \cite{ref2}). As an extension thereof, we asked the following question in 2006 (\cite{ref3} and \cite{ref4}):\\

\textbf{Given any positive integer $N$, can any convex polygonal region $P$ be partitioned into $N$ convex pieces such that all pieces have the same area and the same perimeter?}\\

The convex pieces with same area and perimeter could well have different shapes, different numbers of sides etc. and could be arbitrarily positioned within $P$. To our knowledge, this problem had not been attempted before 2006. We compiled our explorations in \cite{ref5} and this paper is a formal rewrite, mostly.\\

Following \cite{ref4}, we define a \textit{Fair Partition} of a polygon as a partition of it into a finite number of pieces such that every piece has both the same area and same perimeter. Further, if the resulting pieces are all convex, we call it a \textit{Convex Fair Partition}. We may now rephrase our question: Given any positive integer $N$, does every convex polygon allow convex fair partitioning into $N$ pieces?\\

\textit{Remarks:} Any square, rectangle or parallelogram can be convex fair partitioned by $N-1$ parallel lines into $N$ pieces. Any triangle gives $N$ identical triangles if $N$ is a perfect square. Convex fair partitions need not be unique - eg: for $N=4$, a square has infinitely many convex fair partitions. \\

The difficulty with convex fair partitions into $N$ pieces comes from the perimeter constraint; unlike area, we don't know ab initio how much perimeter each piece ought to have; the common value of the perimeter of pieces is an emergent property of the partition.\\

\textbf{The Conjecture:} We tend to believe that every convex polygon allows a convex fair partition into $N$ pieces for any $N$. i.e. every convex polygon can be broken into $N$ convex pieces all of the same area and perimeter. In this paper we discuss only polygonal regions with finite number of sides. But we expect this property to hold for any $2D$ convex region, not only polygonal ones.\\

The conjecture is easily proved for $N=2$ using simple continuity arguments (for example, see \cite{ref5}). Indeed, on the boundary of any convex polygon, there exists \textit{at least one} pair of points \{$P_1$, $P_2$\} such that the line joining them divides the polygon into 2 pieces of equal area and perimeter. Here, we call such a line a \textit{Fair Bisector} of the polygon.\\

\textit{Examples:} (1) for an isosceles triangle with a narrow base, the only fair bisector is the angular bisector of its apex; (2) an equilateral triangle has 3 fair bisectors - its medians; a regular pentagon has 5 fair bisectors. (3) a rectangle has a continuum of infinitely many (any line passing through its center is a fair bisector). If a polygon has a fair bisector between interior points on 2 mutually parallel edges, then it automatically has infinitely many more fair bisectors connecting this pair of parallel edges.\\

\section{Proof of the Conjecture: N=4}

The $N$=2 proof does not readily generalize to $N$=4 and higher powers of 2. Indeed, if one partitions a polygon into 2 pieces of same area and perimeter and recursively does the same to both resulting pieces, one gets, in general, 4 pieces all of the same area, but with the common perimeter of a pair of them different from the common perimeter of the other pair (\cite{ref5} shows a simple numerical example). Below, we show that there is an `enhanced' recursive procedure that does work for $N$=4. \\

\textbf{Definition:} An \textit{Area Bisector} is a straight line that divides a convex polygon into two pieces of equal area (and, in general, different perimeters). From every point on the boundary of a convex polygon, there is an area bisector. A fair bisector is an area bisector giving 2 pieces of equal perimeter. \\

We say a curve $C$ \textit{evolves continuously} with a continuously varying parameter $t$ if (1) every point $P$ on $C$ traces a continuous trajectory as a function of $t$ and (2) the rate of displacement of each point $P$ with change in $t$ varies continuously along $C$. A closed curve stays closed during a continuous evolution but its area and perimeter could both change continuously with $t$.\\

\textbf{The `Augmented' Recursive Scheme for N=4}\\

Divide the `full convex polygon' to be convex fair partitioned into 4 pieces into two equal area pieces by any \textit{area bisector}; call the resulting pieces $A$ and $B$. Consider \textit{a} fair bisector of piece $A$, that divides $A$ into pieces\{$A_1$, $A_2$\} - where both pieces have equal area and equal perimeter. Likewise, consider a fair bisector of $B$ which gives pieces \{$B_1$, $B_2$\} - see figure 1 below; the fair bisectors are shown as dashed lines. In general, at any given initial position of the area bisector, the pair of pieces \{$A_1$, $A_2$\} and the pair \{$B_1$, $B_2$\} have different common perimeters.\\

\begin{figure}[h]
\begin{center}
\hskip 2cm
\includegraphics[width=11.0cm, height=7cm, angle=0]{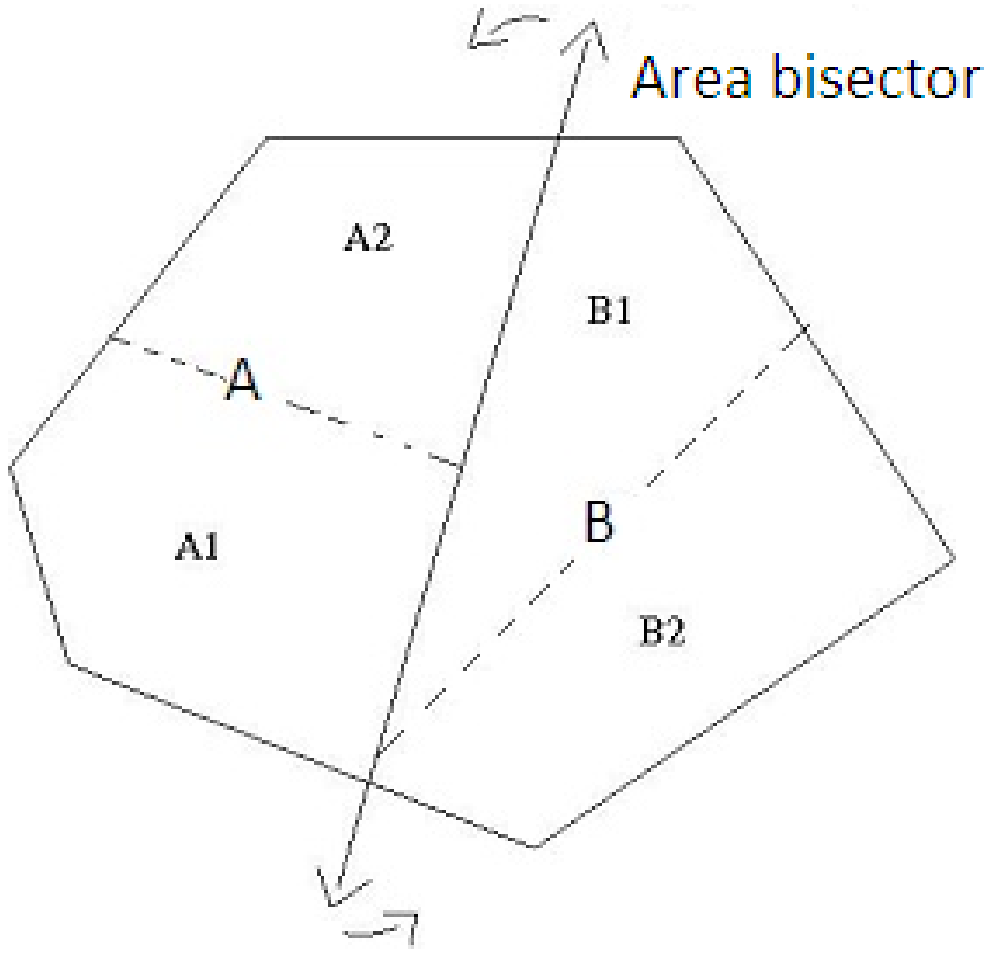}
\caption{\label{} }
\end{center}
\end{figure}
Consider rotating the area bisector of the full polygon from its initial position - by moving both end points of the initial area bisector of the full polygon continuously from their initial positions around the boundary of the polygon so that the 2 pieces $A$ and $B$ continue to have the same area. $A$ and $B$ both change shapes continuously; the fair bisectors of $A$ and $B$ also change (as detailed below). When the area bisector of the full polygon reaches an orientation that is perfectly the reverse of the initial, the piece $A$ would have become the initial $B$ and vice versa. On further rotation of the area bisector, both pieces eventually return to their initial shapes. \\

\textbf{PROPOSITION:} During the full rotation of the area bisector, a state necessarily exists where the 2 pairs of pieces \{$A_1$, $A_2$\} and \{$B_1$, $B_2$\} have the same common perimeter.\\ 

All {$A_1$, $A_2$, $B_1$, $B_2$} - have same area throughout the rotation of the area bisector; so the proposition implies: the augmented recursion achieves a convex fair partition of any convex polygon into 4 pieces.\\

Let the boundary of a convex polygon be linearly parametrized by the length $s$ measured from any boundary point. For every value of $s$, there is a point $P$ on the boundary and also a unique \textit{antipodal point} $P^\prime$ lying  halfway around the boundary from $P$ - and a directed line $PP^\prime$.\\

\textbf{Definition:} The singlevalued function $\alpha(s)$ = difference between area of the part of the polygon lying to the left of (directed) line $PP^\prime$ corresponding to $s$ and area of the part to right of $PP^\prime$.\\

\textbf{Lemma 1:} $\alpha(s)$ is continuous and piece-wise polynomial of degree at most 2 in $s$.\\

\textbf{Proof:} The continuity of $\alpha(s)$ follows from that of the polygon. Let point $P(s)$ vary along an edge, say $E_i$ (which joins vertices $E_{i1}$ and $E_{i2}$ of the polygon) and the corresponding $P^\prime$ move along another edge $E_j$ (joining $E_{j1}$ to $E_{j2}$). Then, $\alpha_l(s)$, the area of the left piece separated by $PP^\prime$ consists of a constant part (the part of the full polygon to the left of the line $E_{i1}$ - $E_{j2}$) and a variable part - the quadrilateral ($E_{i1}$, $P$, $P^\prime$, $E_{j2}$). This varying quadrilateral is, in turn, the union of 2 triangles, ($E_{i1}$, $P$, $E_{j2}$) and ($P$, $P^\prime$, $E_{j2}$). It is easy to see that the area of the former triangle is a linear function of $s$ and area of the latter is at most quadratic in $s$. So $\alpha_l(s)$ is at most quadratic in $s$. Likewise, $\alpha_r(s)$, area of the piece to the right of $PP^\prime$, is quadratic in $s$; and hence also $\alpha(s)$ = $\alpha_l(s)-\alpha_r(s)$. $\diamond$\\

\textit{Note:} Each pair of antipodally located zeros of $\alpha(s)$ corresponds to a fair bisector of the polygon. \\

\textbf{Lemma 2:} For a polygon with finitely many sides, the set of zeros of $\alpha(s)$ is a finite union of points and closed intervals. \\

\textbf{Proof:} By lemma 1, each piece of $\alpha(s)$ is at most of degree 2 in $s$; so in any finite interval of $s$, there will be only finitely many discrete zeros of $\alpha(s)$. There could be continuous intervals of zeros - this happens when instead of the maximum quadratic degree in $s$, the areas of the left and right pieces remain constant and equal for a finite interval of $s$ values (when edges $E_i$ and $E_j$ defined above are parallel). Obviously, the number of such intervals is also finite for a polygon with finitely many sides. $\diamond$ \\

\textbf{Definitions:} A fair bisector is given by a zero of $\alpha(s)$. A \textit{fair range} of a polygon is a connected component of the set of zeros of its $\alpha(s)$. A fair range can either be a point or an interval. Fair ranges occur in antipodal pairs. A fair range is \textit{proper} if $\alpha(s)$ takes different signs on either side of it. A fair bisector is proper if it lies in a proper fair range. Two curves on the same 2D surface have a \textit{proper intersection} if they cut thru each other either at a point or after being coincident in a finite interval; we refer to the length of this interval as the \textit{length of the proper intersection}. Consider a polygon to lie on the $X-Y$ plane and plot values of  its $\alpha(s)$ above the polygon in the $Z$ direction; then, a proper fair range of the polygon is a proper intersection between the polygon and the plot of $\alpha(s).$\\

\textit{Special case:} If $\alpha(s)$ = 0 for \textit{all} values of $s$ (as happens with centrally symmetric polygons - rectangle, regular hexagon, etc..), we consider the full range of $s$ as a single  proper fair range.\\

\textbf{Lemma 3:} The number of proper fair range pairs of any convex polygon is necessarily odd.\\

\textbf{Proof:} From lemma 2, the number of fair ranges is finite. Consider any pair of antipodal points $P$ and $P^\prime$ on the polygon boundary. By definition, $\alpha(s)$ has opposite signs at $P$ and $P^\prime$. This means, if we move continuously from $P$ to $P^\prime$ along the boundary, $\alpha(s)$ will change sign an odd number of times. And since each sign change gives a proper fair range, there is an odd number of them between $P$ and $P^\prime$. If $\alpha(s) = 0$ for all $s$, there is exactly one proper fair range; so there is no exception. $\diamond$\\

We now return to our recursion scheme for $N=4$. Let $\theta$ denote the angle between a reference direction and the rotating area bisector of the full polygon. We focus on one of the resulting pieces, $A = A(\theta)$, which evolves continuously with $\theta$. For piece $A$, consider the function defined above: $\alpha$ = $\alpha_A(s)$ (the subscript `$A$' shows that, in addition to length parameter $s$, $\alpha$ depends on the varying polygon $A$). We consider the proper/improper fair ranges (and corresponding proper/improper fair bisectors) of $A$. Hereon, variable $p$ denotes the common perimeter of pieces \{$A_1$, $A_2$\} cut from $A(\theta)$ by one of $A(\theta)$'s \textit{fair bisectors}. For any value of $\theta$ (a shape of $A$), there could be many values of $p$. On the $\theta-p$ plane, at a particular $\theta$, we can plot fair range of $A(\theta)$ as a continuous interval of $p$ values above that $\theta$.\\

\textbf{Lemma 4:} Let $R$ be a proper fair range of $A$ at $\theta$ = $\theta_0$. There exists a maximal open interval $I$ of $\theta$ values containing $\theta_0$ with the property: for any $b_0 \in R$, there exists a continuous function $b = b(\theta), \theta \in I$ such that $b(\theta_0) = b_0$ and each $b(\theta)$ is an end point of a proper fair bisector of the polygon $A(\theta)$.\\

\textbf{Proof: }Plot $A$ and $\alpha_A(s)$ as given before lemma 3. Proper fair ranges of $A$ are proper intersections between these curves. $A$ and $\alpha_A(s)$ evolve continuously with $\theta$ - with $A$ staying on  $X-Y$ plane. As 2 curves evolve continuously (1) their proper intersections are robust and trace continuous paths (\textit{lengths} of proper intersections - defined above - may vary discontinuously but not their positions); (2) proper intersections do not disappear abruptly but can become tangencies (at a point or along an  interval) and then disappear. Tangencies between $A$ and $\alpha_A(s)$ give improper fair ranges of $A$. Lemma 4 follows. Note: Only tangencies between $A$ and $\alpha_A(s)$ can appear/disappear abruptly with changes in $\theta$. $\diamond$\\

\textbf{Definitions:} At any $\theta$ value, say $\theta_0$, $A$ has at least one proper fair range(from $N=2$ proof); let $R$ be one. As noted above, on $\theta-p$ plane, $R$ gives an interval of values above $\theta_0$ for $p$. By lemma 4, as $R$ changes with $\theta$, its $p$ values give a curve (say) $c$ on the $\theta-p$ plane that can be followed continuously in an interval $I$ of $\theta$ values containing $\theta_0$. We call curves such as $c$, $\gamma$ curves. Moreover as $\theta$ tends to an endpoint, say $\theta_1$, of the interval $I$, $R$ tends to an improper fair range that gives a vertical segment above $\theta_1$ on the $\theta-p$ plane (in special cases, this segment could be just a point); we call such vertical segments, the $\beta$ segments. We define $G$ to be the finite graph given by the union of all $\gamma$  and $\beta$ curves (there are finitely many of each) given by the full $0\leq\theta\leq2\pi$ evolution of $A$.\\

\textbf{Lemma 5:} For each connected component $C$ of $G$, the \textit{parity} of the number of $\gamma$ curve segments above $\theta$ is the same for all $\theta$ values not under a $\beta$ segment.\\

\textbf{Proof:} If, at any $\theta$ value $\theta_1$, $C$ has no $\beta$ segments (for improper fair ranges), all $\gamma$ curves in $C$ continue thru $\theta_1$, so parity of the number of $\gamma$ curves in $C$ cannot change at  $\theta_1$. So, at $\theta = \theta_1$, let there be an improper fair range of $A$ with $\alpha_A(s)$ having zeros at all $s$ values in $[s_1, s_2]$. By definition, $\alpha_A(s)$ has to have values with the same sign (say, positive) for $s$ values in small intervals $[s_1-\epsilon,s_1]$ and $[s_2, s_2+\epsilon]$. Since $\alpha_A(s)$ varies continuously with $\theta$, for $\theta$ close to $\theta_1$, $\alpha_A(s)$ is positive at both $s_1-\epsilon$ and $s_2+\epsilon$, so there are an even number (possibly zero) of sign changes of $\alpha_A(s)$ in $[s_1-\epsilon, s_2+\epsilon]$. Thus, for $\theta$ near $\theta_1$, there are an even number of proper fair ranges within $[s_1-\epsilon, s_2+\epsilon]$. It follows by lemma 4 that on the $\theta - p$ plane, the improper fair range $[s_1, s_2]$ at $\theta = \theta_1$ corresponds to a $\beta$ segment which is met by even numbers of $\gamma$ curves from both sides. The same property obviously holds for all other $\beta$ curves lying above $\theta_1$. It follows that the parity of the number of $\gamma$ curve segments in $C$ lying above $\theta$ is the same as $\theta$ tends to $\theta_1$ from both left and right, and hence stays constant as $\theta$ varies in $[0,2\pi]$. $\diamond$\\

\textbf{Lemma 6:} The graph $G$ contains at least one connected component that spans the full $2\pi$ period of $\theta$ (i.e. whose projection onto the $\theta$ coordinate is $[0,2\pi]$).\\

\textbf{Proof:} If any single fair range of $A$ lasts the full $2\pi$ period, there is nothing to prove. So assume that every proper fair range of $A$ appears at some $\theta$ as an improper fair range and disappears at another $\theta$ as another improper fair range. Consider the set of connected components of $G$. At every $\theta$, $A$ has an odd number of proper fair range pairs (lemma 3); so totally, an odd number of $\gamma$ curves pass over every $\theta$ (obviously, not counting $\gamma$'s that end or begin above that $\theta$). Then, $G$ has at least one connected component $C$ containing an odd number of $\gamma$ segments at that $\theta$. Now, for each connected component of $G$, the parity of the number of $\gamma$ segments is the same for all $\theta$ not under the finitely many $\beta$ segments (lemma 5). So, for component $C$, the number of $\gamma$ segments is odd for all $\theta$ and for no $\theta$ is the number of $\gamma$ segments in $C$ zero (even parity); so $C$ spans $[0,2\pi]$. Indeed, $\beta$ segments in $C$ connect end points of $\gamma$ curves of $C$ ending at same $\theta$, so $C$ contains a continuous curve on $\theta-p$ plane, that spans $[0,2\pi]$. $\diamond$\\ 

\textbf{Lemma 7:} Given a continuous curve $C_1$ defined on the $\theta-p$ plane over $\theta$ such that $C_1$ has period $2\pi$ and spans $\theta = [0, 2\pi]$ (with $\theta = 0$ and $2\pi$ identified). If another curve $C_2$ differs from $C_1$ only by a finite phase $\delta$, then $C_1$ and $C_2$ have at least one common point.\\

\textbf{Proof:} Consider $C_1$ and $C_2$ plotted on the $\theta-p$ plane. Let $C_1$ and $C_2$ reach their (equal) maximum $p$ value, $p_{max}$ and minimum value $p_{min}$ at different values of $\theta$ (else, we already have an intersection). Consider the infinite strip: $0 \leq \theta \leq 2\pi$ on the $\theta-p$ plane. Removal of $(\theta, p)$ points on curve $C_1$ from this strip gives 2 separate semi-infinite strips. Let $C_2$ have the minimum $p$ value $p_{min}$ at point $P_1$ and the value $p_{max}$ at $P_2$ on the $(\theta, p)$ plane. $P_1$ and $P_2$ lie in opposite semi-infinite strips separated by $C_1$. Due to continuity, we can follow $C_2$ from $P_1$ to $P_2$. Curve $C_1$ lies between $p$-values $p_{min}$ and $p_{max}$ and divides the strip $0 \leq \theta \leq 2\pi$ and so, $C_2$ has to cut thru $C_1$ to reach $P_2$ from $P_1$. $\diamond$\\

\textbf{The Final Step:} Let graph $G$, defined before lemma 5, be constructed for both pieces $A$ and $B$. Both graphs obviously have period {$2\pi$} in $\theta$ (the orientation of the area bisector) and an identical sequence of $(\theta, p)$ pairs - since pieces $A$ and $B$ go thru the same sequence of shapes. Indeed, the 2 graphs differ by only a phase $\pi$. By lemma 6, each graph contains a curve spanning the full $[0, 2\pi]$ period. These 2 curves (call them $C_A$ and $C_B$ respectively) differ only by a phase of $\pi$ and by lemma 7, they necessarily intersect within a period. So, there is a value of $\theta$ when perimeters of \{$A_1$, $A_2$\} and those of \{$B_1$, $B_2$\} have a common value. That proves our original proposition and the conjecture for $N=4$. $\diamond$\\

\section{Conclusion}

A sophisticated proof of our conjecture for $N=3$ has been published (\cite{ref6}). Some arguments for $N=3$ are also in \cite{ref5}. In \textit{appendix} below, we show that our $N=4$ proof generalizes to all powers of 2. We suspect, if the conjecture is proved for any $N$ (say, for 3, as in \cite{ref6}), our arguments could yield a proof for $2N$. We also believe our $N=4$ proof for polygonal regions can be extended to any 2D convex region.\\

This is known: if we need only to \textit{fair partition} polygons into $N$ pieces (ie the pieces have same area and perimeter but don't have to be convex), it is possible for any (even non-convex) polygon and any $N$ (\cite{ref5}). By simple examples, we see that the fair partition which \textit{minimizes} the total perimeter of pieces is \textit{not} \textit{necessarily} a convex fair partition (\cite{ref5}). Such `optimal' fair partitioning seems an unexplored area. \\

If a smart counter example is found disproving our conjecture for some $N$, we could ask: how to \textit{decide} if a convex fair partition exists for a given convex polygon and given $N$? Finally, what about higher dimensional analogs of this problem?\\

\textbf{Appendix: $N$ = 8 and higher powers of 2}\\

\textbf{Claim on $N$ = 8:} Divide the polygon to be 8-partitioned, say $P$, into 2 continuously changing pieces, $X$ and $Y$, by a rotating area bisector (whose orientation is given by angle $0\leq\phi<2\pi$). If we 4-partition both $X$ and $Y$ at every $\phi$, then there is some value of $\phi$ when the 2 sets of 4 (equal area) pieces from both $X$ and $Y$ have the same common perimeter, thus proving the conjecture for $N=8$. \\

\textbf{Proof (sketch):} We begin by stating (proof omitted) a special case of lemma 7: \\

\textbf{Lemma 8:} On a plane, say $\theta-p$, given a continuous curve $C_1$ defined over $\theta$ such that $C_1$ has period $2\pi$ and spans $\theta = [0, 2\pi]$ (with $\theta = 0$ and $2\pi$ identified). If another curve $C_2$ differs from $C_1$ only by a phase $\pi$, then plots of $C_1$ and $C_2$ have an \textit{odd number of pairs} of proper intersections in $0\leq\theta<2\pi$ - with both points in a proper intersection pair having same $p$ values and their $\theta$s differing by $\pi$.\\

Consider piece $X(\phi)$ cut from $P$ by a rotating area bisector. At each $\phi$, $X$ is 4-partitioned as in $N=4$: $X$ is divided into $A$ and $B$ by $X$'s area bisector with orientation $\theta$. $p$ is a common perimeter of pieces \{$A_1$, $A_2$\} cut from $A(\theta)$ by a fair bisector of $A(\theta)$. On the $\theta-p$ plane is the periodic curve $C_A(\theta)$ (defined in 'final step' above) tracing the evolution of $A(\theta)$. Curve $C_B(\theta)$ (a copy of $C_1$ with phase difference $\pi$) similarly follows $B(\theta)$. Each intersection between $C_A$ and $C_B$ gives a 4-partition of piece $X$. \\

Now, in our $N=8$ method, at each value of $\phi$, $\theta$ makes a complete $2\pi$ rotation and curves $C_A$ and $C_B$ (both defined for piece $X$ at that $\phi$), follow the change of $\theta$. Now, if $\phi$ also varies continuously, both $C_A$ and $C_B$ evolve continuouly and proper intersections between $C_A$ and $C_B$ trace continuous trajectories. By lemma 8, at every $\phi$, the number of proper intersection pairs of $C_A(\theta)$ and $C_B(\theta)$ is odd. This, with lemma 6, implies that a continuous, periodic curve, say $C_X(\phi)$, can be formed from paths traced by proper intersections of $C_A$ and $C_B$ as $\phi$ changes. A similar curve $C_Y(\phi)$ exists for piece $Y$. Obviously, $C_Y$ is identical to $C_X$ except a phase of $\pi$. By lemma 7, $C_X$ and $C_Y$ intersect. Since every point of $C_X$ ($C_Y$) gives a 4-partition of $X$ ($Y$), intersections of $C_X$ and $C_Y$ give 8-partitions of $P$.$ \diamond$\\

Recursive applications of the above arguments prove the conjecture for all powers of 2.\\

\textbf{Acknowledgements:}
Discussions with John Rekesh, Pinaki Majumdar, Arun Sivaramakrishnan, Varun Vikas, Bhalchandra Thatte and researchers at DAIICT Gandhinagar were a great help. Thanks to Kingshook Biswas for his guidance and advice. \\


\begin{thebibliography}{99}

\bibitem{ref1} Akiyama, Kaneko, Kano, Nakamura, Rivera-Campo, Tokunaga, and Urrutia. `Radial perfect partitions of convex sets in the plane'. In Japan Conf. Discrete Comput. Geom., pages 1-13, 1998. 

\bibitem{ref2}  Jin Akiyama, Gisaku Nakamura, Eduardo Rivera-Campo, and Jorge Urrutia. 
`Perfect divisions of a cake' In Proc. Canad. Conf. Comput. Geom., pages 114-115, 1998.

\bibitem{ref3} http://nandacumar.blogspot.com/2006/09/cutting-shapes.html (September 2006)

\bibitem{ref4} http://maven.smith.edu/~orourke/TOPP/P67.html (July 2007)

\bibitem{ref5}  R. Nandakumar and N. Ramana Rao. http://arxiv.org/abs/0812.2241v2 (December 2008).

\bibitem{ref6}  I. Barany, P. Blagojevic and A. Szucs.  `Equipartitioning by a Convex 3-fan'.
Advances in Mathematics, Volume 223, Issue 2, January 2010, pages 579-593.

\end{thebibliography}
\end{document}